\documentclass{article}
\usepackage[all]{xy}
\usepackage{graphicx}
\usepackage{amssymb}
\usepackage{amsthm}
\renewcommand{\qedsymbol}{\rule{5pt}{5pt}}
\usepackage{mathtools}
\usepackage{mathrsfs}
\usepackage{abstract}
\title{\bf Equivalence of Weak and Strong Uniform Glivenko-Cantelli Classes}
\author{Shane R. Crowe}
\begin{document}
\maketitle

\begin{abstract}
Strong uniform Glivenko-Cantelli classes are weak uniform Glivenko-Cantelli, but to date results that prove the converse require additional conditions on the function class. This paper shows that no such extra requirements are necessary.
\end{abstract}

\section*{Introduction}

Let ($S, \Sigma)$ be a measurable space, $\mathcal{P}$ the class of all probability measures on it, and $H$ a class of real valued integrable functions on $S$, assumed throughout to contain at least one non-constant function in order to avoid trivialities. Then we say that $H$ is a {\it weak uniform Glivenko-Cantelli} class 
if for any $\epsilon > 0$ we have

\[ \lim_{n \rightarrow \infty} \sup_{P \in \mathcal{P}} P^n\{||P_n - P||_H^* > \epsilon\} = 0 \]

\noindent
and a {\it strong uniform Glivenko-Cantelli} class if for any $\epsilon > 0$ we have

\[ \lim_{n \rightarrow \infty} \sup_{P \in \mathcal{P}} P^{\infty}\{\sup_{k \geq n} ||P_k - P||_H^* > \epsilon\} = 0 \]

\noindent
where asterisks denote measurable cover functions, $k, n \in \mathbb{N}$, $(S^n,\Sigma^n, P^n)$ is the usual product space, $P_n$ is the empirical probability on $n$ i.i.d samples under $P$, and $||.||_H$ is given by

\[ ||F||_H = \sup_{f \in H} |F(f)| \]

\noindent
for any real valued functional $F$.

\vspace{5.0mm}

\noindent
If we drop the requirement of uniformity of convergence over $\mathcal{P}$, so we ask that convergence to zero in both definitions holds only for all elements of $\mathcal{P}$ individually, then we call $H$ (respectively) a weak/strong {\it universal} Glivenko-Cantelli class, and it is known $[2]$ that in fact these notions are equivalent,
i.e. a weak universal Glivenko-Cantelli class is strong universal Glivenko-Cantelli and vice versa. A strong uniform Glivenko-Cantelli class $H$ is of course weak uniform Glivenko-Cantelli, but the converse is only known to hold if additional conditions are placed on $H$, for example that it be image admissible Suslin $[1],
[2]$. This paper will show that in fact the converse holds without any such special requirements on $H$.

\section*{Background Material}

{\bf Proposition:} Suppose $H$ is a weak uniform Glivenko-Cantelli class. Then $0 <  M < \infty$ where

\[ M = \sup_{f \in H}[\sup f - \inf f]  \]

\noindent
This is proved in Proposition 10.2 in $[1]$. It follows immediately that for a weak uniform Glivenko-Cantelli class $H$ that

\[ |f(x) - f(y)| \leq M \]

\noindent
for any $f \in H$ and any $x, y \in S$, and that for any $n$ and any $P \in \mathcal{P}$ that

\[ ||P_n - P||_H = \sup_{f \in H}|P_n(f-\inf f) -  P(f - \inf f)|  \leq 2M \]

\vspace{5.0mm}

\noindent
{\bf Lemma:}  If $H$ is a weak uniform Glivenko-Cantelli class then 

\[ \lim_{n \rightarrow \infty} \sup_{P \in \mathcal{P}} E[||P_n - P||_H^*] = 0 \]

\noindent
{\bf Proof:} For any $n$ and $P \in \mathcal{P}$ we have that $||P_n - P||_H^* \leq 2M$ with $P^n$ probability $1$. Then for any $n$ and any $\epsilon > 0$

\begin{flalign*}
\sup_{P \in \mathcal{P}} E[||P_n - P||_H^*] &= \sup_{P \in \mathcal{P}}[ \int_{\{||P_n - P||_H^* > \epsilon\}} ||P_n - P||_H^* dP^n &\\
&+ \int_{\{||P_n - P||_H^* \leq \epsilon\}} ||P_n - P||_H^* dP^n]&\\
&\leq  2M \sup_{P \in \mathcal{P}} P^n\{||P_n - P||_H^* > \epsilon\} + \epsilon  &\\
\end{flalign*}

\noindent
and the result now follows from the assumption that $H$ is weak uniform Glivenko-Cantelli. $\qedsymbol$

\vspace{5.0mm}

\noindent
For any $n$ and $P \in \mathcal{P}$ the function $||P_n - P||_H$ satisfies McDiarmid's bounded difference property with constants $Mn^{-1}$, which is to say that if $s = (s_1,...,s_n) \in S^n$ is any point, and for any $1 \leq i \leq n$ we exchange $s_i$ with any $y \in S$, calling this point $s'$ say, then

\[ | \,||P_s - P||_H - ||P_{s'} - P||_H \, | \leq \sup_{f \in H} |P_sf - P_{s'}f| \]

\noindent
(where $P_s, P_{s'}$ denote the empirical probabilities at the subscripted points),

\[  =  \sup_{f \in H} |n^{-1}(f(s_i) - f(y))|  \leq Mn^{-1} \]

\noindent
For any $n$ and $P \in \mathcal{P}$ we have (after a possible null set modification) that $||P_n - P||_H^* \leq 2M$, and then following the definition in $[3]$, the function $||P_n - P||_H^*$ is weakly difference-bounded by $(4M, Mn^{-1}, 0)$ because if
$s = (s_1,...,s_n) \in S^n, y \in S$ are any points and $s'$ is $s$ with $s_i$ replaced by $y$ again then

\[ | \, ||P_s - P||_H^* - ||P_{s'} - P||_H^* \,  | = | \, [ \, ||P_{\bullet} - P||_H^* \circ \pi](s,y) - [ \, ||P_{\bullet} - P||_H^* \circ \pi \circ \sigma_i](s,y) \, | \]

\noindent
where $\pi(s_1,...,s_n,y) = (s_1,...,s_n)$ is a projection, and $\sigma_i$ the permutation that swaps the $i$-th and $(n+1)$-th entries in any elt of $S^{n+1}$. Then since projections and permutations on i.i.d. spaces are perfect functions we have almost surely that

\[ | \, ||P_s - P||_H^* - ||P_{s'} - P||_H^* \, | \leq | \, ||P_s - P||_H - ||P_{s'} - P||_H \, |^*   \]

\noindent
and the RHS is again bounded by $Mn^{-1}$ a.s.

\section*{Main Result}

\vspace{5.0mm}

\noindent
{\bf Theorem:} A weak uniform Glivenko-Cantelli class is strong uniform Glivenko-Cantelli.

\vspace{5.0mm}

\noindent
{\bf Proof:}  Let $H$ be a weak uniform Glivenko-Cantelli class. As discussed above, for any $n$ and $P \in \mathcal{P}$ we have that $||P_n - P||_H^*$ is weakly difference-bounded by $(4M, Mn^{-1},0)$ and so by Corollary 4.7 in $[3]$ for any $\epsilon > 0$ we can say that

\[ P^n\{||P_n - P||_H^* - E[||P_n - P||_H^*] > \epsilon/2\} \leq  2e^{-3n\epsilon^2/8M(15M + \epsilon) } \]

\noindent
Therefore

\[  P^n\{||P_n - P||_H^* - \sup_{P \in \mathcal{P}} E[||P_n - P||_H^*] > \epsilon/2\} \leq 2e^{-3n\epsilon^2/8M(15M + \epsilon) } \]

\noindent
The Lemma from earlier tells us that $\sup_{P \in \mathcal{P}} E[||P_n - P||_H^*] \rightarrow 0$, so there is an $n_0 \in \mathbb{N}$ such that $n > n_0$ implies that $\sup_P E[||P_n - P||_H^*] < \epsilon/2$ and thus $n > n_0$ implies that

\[  P^n\{||P_n - P||_H^*  > \epsilon\} \leq 2e^{-3n\epsilon^2/8M(15M + \epsilon) } \]

\noindent
for all $P \in \mathcal{P}$. To show that $H$ is strong uniform Glivenko-Cantelli we need to prove that the sequence

\[ \sup_{P \in \mathcal{P}} P^{\infty}\{\sup_{k \geq n}||P_k - P||_H^* > \epsilon\} \]

\noindent
converges to zero, for any $\epsilon > 0$. Note that for any $n$

\begin{flalign*}
\sup_{P \in \mathcal{P}} P^{\infty}\{\sup_{k \geq n}||P_k - P||_H^* > \epsilon\} &= \sup_{P \in \mathcal{P}} P^{\infty}[\bigcup_{k \geq n} \{||P_k - P||_H^* > \epsilon\}] &\\
&\leq \sup_{P \in \mathcal{P}} \sum_{k \geq n} P^k\{||P_k - P||_H^* > \epsilon\} &\\
\end{flalign*}

\noindent
Thus for $n > n_0$ we have that 

\begin{flalign*}
\sup_{P \in \mathcal{P}} P^{\infty}\{\sup_{k \geq n}||P_k - P||_H^* > \epsilon\} &\leq 2\sum_{k \geq n} e^{-3k\epsilon^2/8M(15M + \epsilon) }&\\
& = 2Le^{-3n\epsilon^2/8M(15M + \epsilon) } &\\
\end{flalign*}

\noindent
where $L < \infty$ is the sum of the geometric series $\sum_{j=0}^{\infty} e^{-3j\epsilon^2/8M(15M + \epsilon) }$, and the RHS above converges to zero, as required. $\qedsymbol$

\section*{References}

\noindent
[1] R.M. Dudley, Uniform Central Limit Theorems (2014). {\it Cambridge University Press}.

\vspace{5.0mm}

\noindent
[2] R.M. Dudley, E. Gin\'e, and J. Zinn, Uniform and Universal Glivenko-Cantelli Classes (1991). {\it Journal of Theoretical Probability} {\bf 4}, 485–510.

\vspace{5.0mm}

\noindent
[3] S. Kutin, Extensions to McDiarmid’s inequality when differences are bounded with high probability (2002). Technical Report TR-2002-04, {\it Department of Computer Science, University of Chicago}.

\end{document}